\documentclass[a4paper,11pt]{article}
\newcommand{\beq}{\begin{equation}}
\newcommand{\eeq}{\end{equation}}

\usepackage{amsthm,amssymb,amsmath}
\usepackage{graphicx}
\usepackage{wasysym}
\usepackage{epstopdf}
\usepackage{mathdots}
\usepackage[lmargin=.75in,rmargin=.75in,tmargin=.75in,bmargin=.75in]{geometry}
\usepackage{float}
\usepackage{hyperref}
\usepackage{subfigure}

\theoremstyle{remark}

\theoremstyle{definition}

\begin{document}

\title{GPU-Acceleration of Parallel Unconditionally Stable\\
Group Explicit Finite Difference Method}
\author {K. Parand $^{\mbox{\footnotesize }}$\thanks{E-mail
address: \texttt{k\_parand@sbu.ac.ir}, Corresponding author, (K. Parand)}\hspace{1mm}
, \normalsize{Saeed Zafarvahedian$^{\mbox{\footnotesize
}}$\thanks{E-mail address: \texttt{zafarvahedian@iran.ir; S.Zafarvahedian@Mail.sbu.ac.ir},
(S. Zafarvahedian)}\hspace{1mm}}
, \normalsize{Sayyed A. Hossayni$^{\mbox{\footnotesize
}}$\thanks{E-mail address: \texttt{hossayni@iran.ir},
(Sayyed A. Hossayni)}\hspace{1mm}}
\vspace{-2mm}\\\\
\footnotesize{\em $^{\mbox{\footnotesize a}}$Department$\!$ of$\!$
Computer Sciences, Faculty$\!$ of$\!$ Mathematical Sciences$\!$,}\vspace{-2mm}\\
\footnotesize{\em Shahid Beheshti University, Evin,Tehran 19839,Iran}\\
} \maketitle
\vspace{-.9cm}\noindent\linethickness{.5mm}{\line(1,0){475}}
\maketitle
\begin{abstract}
Graphics Processing Units (GPUs) are high performance co-processors originally intended to improve the use and quality of computer graphics applications.
Once, researchers and practitioners noticed the potential of using GPU for general purposes, GPUs applications have been extended from graphics applications to other fields.
The main objective of this paper is to evaluate the impact of using GPU in solution of the transient diffusion type equation by parallel and stable group explicit finite difference method and encourage the researchers in this field to immigrate from implementing their algorithms in CPU to the GPU emerging world.
For comparing them, we implemented the method in both GPU and CPU (multi-core) programming context.
Moreover, we proposed an optimal synchronization arrangement for the implementation pseudo-code.
Also, the interrelation of GPU parallel programming and initializing the algorithm variables were discussed, taking advantage of numerical experiences.
The GPU-approach results are faster than those obtained from a much expensive parallel 8-thread CPU-based programming.
The GPU used in this paper, is an ordinary old laptop GPU (GT 335M, launched at 2010) and is accessible for everyone and the newer generations of GPU (as discussed in paper) have even more performance priority over the similar-price GPUs.
Then, the results are expected to encourage the entire research society to take advantage of GPUs and improve the time efficiency of their studies.
\vspace{.3in}
\end{abstract}
\section{Introduction}
\subsection{Graphics processing unit}
To attain a high-level graphical quality with adequate performance, computer graphical applications needed to a very fast computation of large amounts of data.
For this sake, some dedicated high-performance multi-processors, named Graphics Processing Units (GPU) were widely utilized.
A GPU is a Single Instruction Multiple Data (SIMD), parallel architecture, in which one instruction can be performed with multiple data in a pipeline, simultaneously\cite{Heimlich.2011}.
As scientists and engineers faced with such high computing potential in GPUs, they started to use it for general purposeful projects\cite{Kyoung.2004,Hsien.2005}.
Not only for their great speed up, but also for their low cost platform, developing GPU-based applications, in several cases, were preferred to the CPU-based applications.\\
Gradually, programmable GPUs became processors with several cores and threads, excellent computational speed and high-level memory bandwidth.
Therefore, aggregation of features like general purpose super computability, high parallelism, very high memory bandwidth, low cost and compact size, converted a GPU-equipped Laptop to a successful substitute for even a massive parallelism CPU-based system\cite{Mielikainen.2011}.
CUDA is an add-on for C programming language by which a programmer can directly run instructions on GPUs. CUDA can, also, help scientists and engineers to access GPU as they accessed CPU before, by simple programming codes.
Nowadays, GPUs are able to execute various science and engineering algorithms, up to more than 100 times faster than similar-cost CPUs \cite{Abdellah2015} (when the CPU code is single-threaded) and up.
Nowadays, almost every researcher has a programmable GPU inside his computer and may not be advised about this unused power, which can be accessed by simple and similar programming codes to what they did before.
In this paper, we aim to accelerate "the fully explicit unconditionally stable parallel finite difference method", using an ordinary laptop GPU (GT-335M) released in 2007.
So, the results are expected to encourage the researchers in this field to use GPU for heavy computations, required in their studies.
However, the advanced potential of GPUs increases the number of scientists and engineers who tend to apply this as an affordable and efficient platform, day by day.\\
Oh\cite{Kyoung.2004} has used GPU for a faster artificial neural network. It is used to implement the matrix multiplication of a neural network to enhance the time performance of a text detection system.
The writers of \cite{Georgii.2005} have presented and analyzed different implementations of mass-spring systems for interactive simulation of deformable surfaces on GPUs.
Ho\cite{Ho.2006} presents practical and real-time rendering methods for  illumination adjustable image (IAI) based on spherical Gaussian kernel functions.
Heimlich et al.\cite{Heimlich.2011} evaluated the impact of using GPU in two computational expensive problems in nuclear area:
a neutron transport simulation by Monte Carlo method and Poisson heat equation resolution for a square thin plate.
Liu et al. \cite{liu2014solving} presented a parallel GPU solution for the Caputo fractional reaction-diffusion equation in one spatial dimension with explicit finite difference approximation.
\subsection{Group explicit finite difference method}
Most of the transport phenomena in the scientific computation have usually transportation of special property with the diffusion\cite{Davam.2007}.
In this paper, we consider the numerical solution of diffusion type transient equation$(U_t=\nabla^2U)$.
Although, practically, the transportation with diffusion may be mixed with other transportation mechanisms, every solution which is capable of solving pure diffusion equations may directly be used for complex transport equations, as well.
Also, the diffusion step may, separately, be solved by operator splitting method\cite{Davam.2007}.\\
Finite difference is a popular numerical method for solution of the transport equations.
In this method, the space is discretized using structured grids and the time axis is managed using explicit or implicit techniques\cite{Davam.2007}.
Numerical models, based on explicit time discretization, have some benefits for parallelization.
They make the parallelization process easier.
Also, taking their advantage, we can reach higher parallelization performances in comparison with similar methods. However, due to the stability need, they are under high restriction in time step.
Using other parallelization approaches which depend on implicit techniques and do not have time step restriction, we should solve a linear or nonlinear algebraic equations system, in each step, and obviously it is not straightforward to implement them on parallel computers, at all.
Moreover, since ordinary equation systems are often solved by iterative methods, the necessary iteration number to reach an acceptable accuracy may be large; especially for long time increments and small sized meshes.\\
In the last decade, a so-called alternating explicit-implicit schemes were studied by a large number of researchers. This scheme is based on the domain decomposition concept and is a composition of implicit and explicit schemes\cite{Davam.2007}.
There, a physical domain is divided into several sub-domains.
In each step, the problem was solved explicitly in some sub-domains and implicitly in others.
Then, in the following step, the scheme is reversed (Those which were used to be solved implicitly will be solved explicitly and vice versa).
One of the advantages of explicit-implicit schemes is parallelization.
Although some members of this family are unconditionally stable, some others suffer some limitations like time step restriction.
However, these restrictions are much lower than those of a fully explicit method.
Notwithstanding good stability and intrinsic parallel nature of alternating explicit-implicit schemes, these methods have, still, implicit phase which forces implicit solution of problems over implicit sub-domains and this leads to the solution of a system of linear or non-linear algebraic equations at each time level.\\
Although considering the local nature of such equations, there is no problem for parallel implementation of the related system, in general, it is preferable to avoid solving algebraic equation systems.\cite{Davam.2007}
Moreover, the extension of these schemes to higher dimensions (rather than one dimension) makes some computational complexities.\\
The need to parallel unconditionally stable and fully explicit methods led Tavakoli et al. \cite{Davam.2007} to develop such scheme.
\subsection{Paper achievements}
In this paper, we implement Tavakoli's \cite{Davam.2007} parallel and stable explicit finite difference method for the solution of one-dimensional diffusion equation by a multi-thread GPU programming code.
To accomplish that, first of all, we discuss on some techniques for parallel implementation of the method and discuss on the related issues such as detecting the optimized location for synchronizations and etc.
Then, before running the algorithm, we investigate how to set the values of the algorithm variables.
Thereafter, we provide comparisons with CPU (single and multi-thread) approaches both of which were accomplished by an ordinary Laptop.
For this purpose, the CPU multi-thread and single-thread programming codes were also implemented, alongside, the GPU code, itself.
\subsection{The paper hardware specification}\label{hardware.sec}
The GPU used in this paper is GeForce GT-335M, the third generation of the CUDA enabled NVIDIA GPUs.
GT-335M has 9 blocks and 8 stream processors per each block (totally 72 stream processors) and is capable of demonstrating the result of 233 GFLOPS \cite{geForcePerf}.
The GT-335M frame buffer memory interface is 128 bits wide, composed by GDDR3 memory controllers.
The controllers coalesce a much greater variety of memory access patterns, improving the efficiency as well as peak performance\cite{Heimlich.2011}.
Also, the memory controller is configured with 1 GB of memory.
The eager reader to GPU structure and more information about GPU components, is referred to \cite{Heimlich.2011}
The CPU, used in this study, is Intel Core i7-720QM, one of the Intel Core i7-700 mobile processors series and is capable of demonstrating the result of 25-45 GFLOPS \cite{intelPerf}.
It has 4 cores and 8 threads and is composed by DDR3 memory controllers.
Note that since the audience of this study are the broad society of physics and chemistry who are not expert in computer science and in complicated programming, we avoided using SSE (Streaming SIMD Extensions) \cite{kusswurm2014streaming} instructions of the CPU \cite{intelSSE} and also do not recommend using SFU (Special Function Unit) of GPU \cite{lindholm2008nvidia}, each of which can increase the performance of the code in CPU and GPU, but also increase the coding complexity.
In 2010, the price of GT-335M was about 250\$ while the price of Intel Core i7-720QM was about 380\$.
All of the codes have been run by both processors (written normally in C language) and the results show GPU excellence in spite of its cheaper price in comparison with CPU.
\section{The parallel scheme for one-dimensional diffusion equation}
\begin{figure}
\centering
\subfigure[The computational molecule for the Saul’yev formulae A and B.]{
\centering
\includegraphics[scale=0.4]{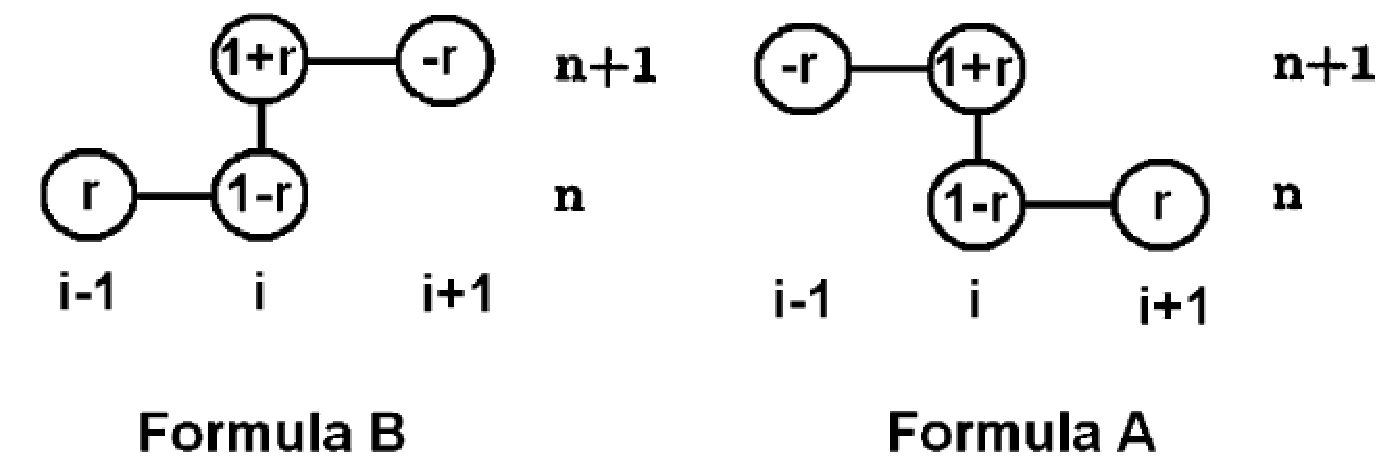}
\label{AB}
}
\subfigure[The computational molecule for nodes $x_{2ik+k-1}$ ((a)) and nodes $x_{2ik+k}$ ((b)).]{
\centering
\includegraphics[scale=0.20]{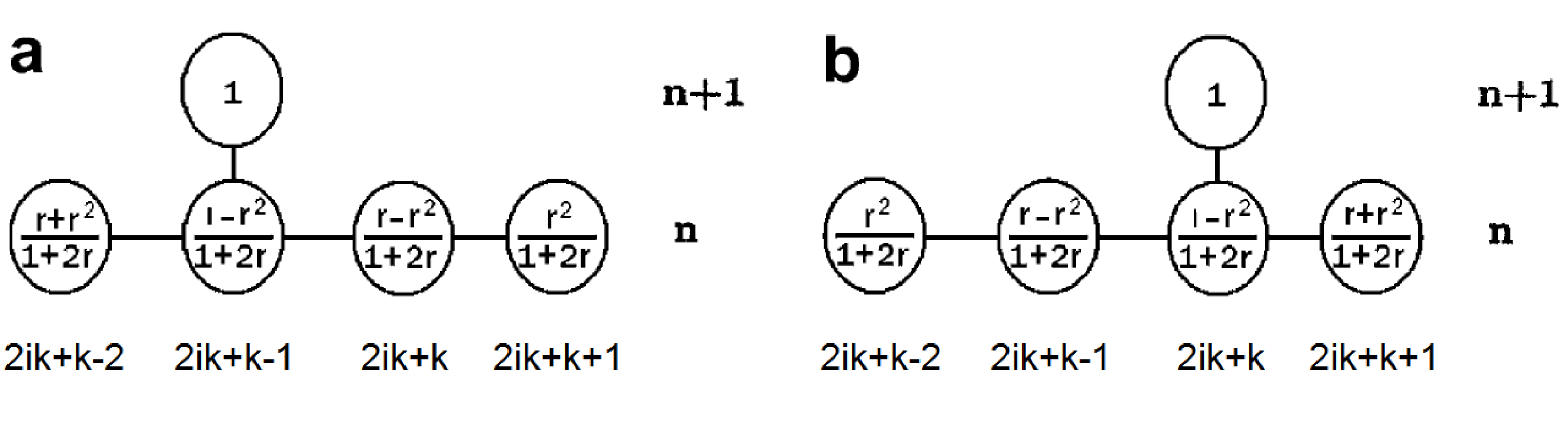}
\label{aabb}
}
\subfigure[The diagram of the parallel implementation of asymmetric Saul’yev schemes for two sub-domains.]{
\centering
\includegraphics[scale=0.21]{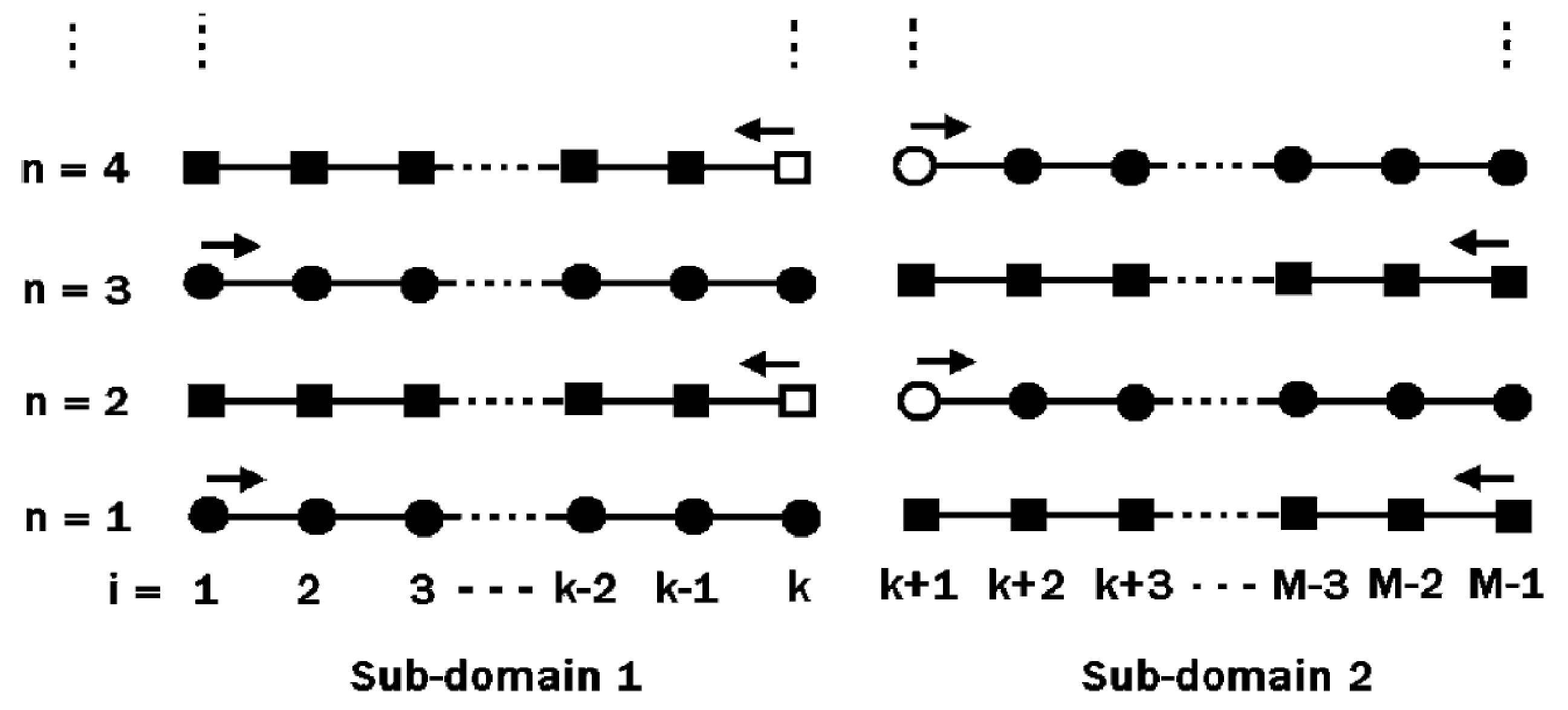}
\label{f_par_tav}
}
\subfigure[The diagram of the parallel implementation of asymmetric Saul’yev schemes.]{
\centering
\includegraphics[scale=0.8]{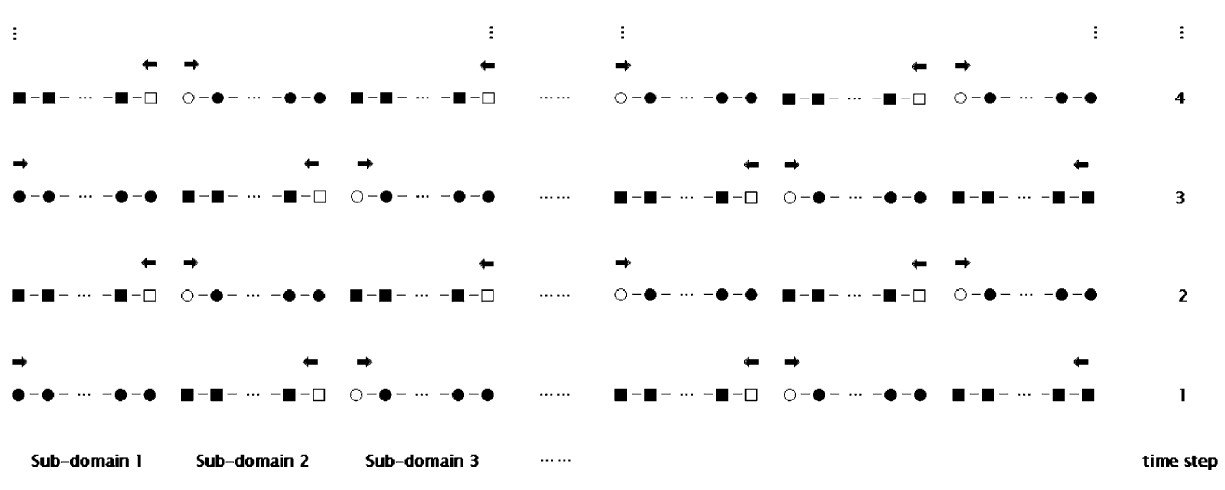}
\label{par_tav}
}
\caption{A schematic view for parallel and stable group explicit finite difference method}
\end{figure}

This section addresses the solution of the following problem
\begin{align*}
\begin{cases}
\frac{\partial u}{\partial t}=\frac{\partial^2u}{\partial x^2},&(x,t)\in \Omega\cr
u(x,0)=\phi(x),&0\leq x\leq 1\cr
u(0,t)=u(1,t)=\beta&0\leq t\leq T
\end{cases}
\end{align*}
where $\Omega= [0,1]\times [0,t]$, $\phi(x)$ is initial condition and $\beta$ is prescribed constant boundary condition.
As we mentioned, we are going to implement Tavakoli's\cite{Davam.2007} method by means of GPU and CPU multi-thread programming.
So, we should, first, introduce his method.\\
As it can be seen in \cite{Davam.2007}, the domain $[0,1]\times[0,T]$ is divided into $M\times N$ finite difference scheme mesh nodes with spatial step size $h = \frac{1}{M}$ in $x$ direction and the time step size $\tau = \frac{T}{N}$ respectively.
Grid points $(x_i, t_n)$ are given by
\begin{align*}
&x_i = ih,~~~i = 0,1,2,. . . ,M,&&t_i = n\tau,~~~n =0, 1, 2, . . . ,N,
\end{align*}
in which $M$ and $N$ are positive integers.
$u_i^n$ is used to denote the finite difference approximations of $u\left(i h, n \tau\right)$.
Using Saul'yev's two schemes for the second derivative $\left(\frac{\partial^2u}{\partial x^2}\right)$, Tavakoli\cite{Davam.2007} proposed two relations for walking on the finite difference scheme mesh nodes.
The first relation which is extracted from Saul'yev scheme A\cite{Davam.2007} is
\begin{equation}\label{aa}
u_i^{n+1}=\frac{1}{1+r}\left[r u_{i-1}^{n+1}+(1-r)u_i^n+r u_{i+1}^n\right]
\end{equation}
for $i = 1,2,\cdots,M-1$ where $r=\frac{\tau}{h^2}$.
He has noted that the Saul’yev scheme A is explicit if the calculation begins at the left boundary, and moves to the right, so that, only the single value $u_i^{n+1}$ is unknown.
The second equation, extracted from Saul'yev scheme B\cite{Davam.2007}, is
\begin{equation}\label{bb}
u_i^{n+1}=\frac{1}{1+r}\left[r u_{i+1}^{n+1}+(1-r)u_i^n+r u_{i-1}^n\right]
\end{equation}
for $i = M-1,M-2,\cdots,2,1$.
The Saul’yev scheme B has a similar explicit nature if the calculation proceeds to the left from the right hand boundary.
The computational molecule of this formula is shown in figure \ref{AB}.
Application of the Saul’yev formula A and B alternatively provides better accuracy and error distribution.\\
For parallel implementation of Saul’yev's scheme, Tavakoli\cite{Davam.2007} used the concept of domain decomposition and proposed a preliminary scheme, illustrated in the figure \ref{f_par_tav}.
He, also, proposed a complementary scheme which is development of the preliminary one, to more than two sub-domains.
Description of the complementary scheme would illustrate both the former and the latter schemes for the reader.
So, we neglect much description of the former and, only, describe the latter.\\
For convenience, let $M$ to be an even multiple of $k$ (e.g. $M=2\*e\*k$, for arbitrary $k$); so, it decomposes the domain $\Omega$ into, respectively, odd and even sub-domains:
\begin{align*}
&[x_{2ik},x_{2ik+k-1}],\\
&[x_{2ik+k},x_{2(i+1)k-1}],\\
&i=0,1,\cdots,e-1.
\end{align*}
If, as depicted in \ref{par_tav}, the formula A is used for the odd sub-domains $[x_{2ik},x_{2ik+k-1}]$ and the calculation begins at the left wall of the sub-domain and moves to the right (and, in a similar manner, the formula B is used for the even sub-domains $[x_{2ik+k},x_{2(i+1)k-1}]$ and the calculations proceed to the left from the sub-domain right wall), then calculating each sub-domain is independent from the other.
For the next time-step, the formula B is used for the odd sub-domains and the formula A is used for the even sub-domains and the calculation is started from the points $x_{2ik+k-1}$ and $x_{2ik+k}$ for odd and even sub-domains, respectively.\\
Moreover, Tavakoli\cite{Davam.2007} has calculated two formulas for confluence points.
At the nodes $(x_{2ik+k-1}, t_{n+1})$, he has approximated the solution by
\begin{align}\label{cc}
u_{2ik+k-1}^{n+1}=\frac{1}{3(1+2r)}
\left[(r+r^2)u_{2ik+k-2}^n+(1-r^2)u_{2ik+k-1}^n+(r-r^2)u_{2ik+k}^n+r^2u_{2ik+k+1}^n\right]
\end{align}
while at nodes $(x_{2ik+k}, t_{n+1})$, he has approximated it by
\begin{align}\label{dd}
u_{2ik+k}^{n+1}=\frac{1}{3(1+2r)}
\left[r^2u_{2ik+k-2}^n+(r-r^2)u_{2ik+k-1}^n+(1-r^2)u_{2ik+k}^n+(r+r^2)u_{2ik+k+1}^n\right]
\end{align}
which can be diagrammatically represented by the 2 computational molecules given in figure \ref{aabb}.
After updating the nodes $x_{2ik+k-1}$ by Eq. (\ref{cc}) and the nodes $x_{2ik+k}$ by Eq. (\ref{dd}) the calculation of the remained nodes will continue by Eq. (\ref{bb}) for odd subdomains and by Eq. (\ref{aa}) for even sub-domains.
These sequence computation is repeated for each pair of consecutive time steps.
The diagram of this method is displayed in figure \ref{par_tav}.
In Tavakoli's article\cite{Davam.2007}, the symbol $\CIRCLE$ has been applied to denote the scheme Eq. (\ref{aa}), $\blacksquare$ for the the scheme Eq. (\ref{bb}), $\Box$ for the scheme Eq. (\ref{cc}) and $\Circle$ for the scheme Eq. (\ref{dd}).
\section{Implementation of the parallel scheme}
\begin{figure}
\centering
\includegraphics[scale=0.30]{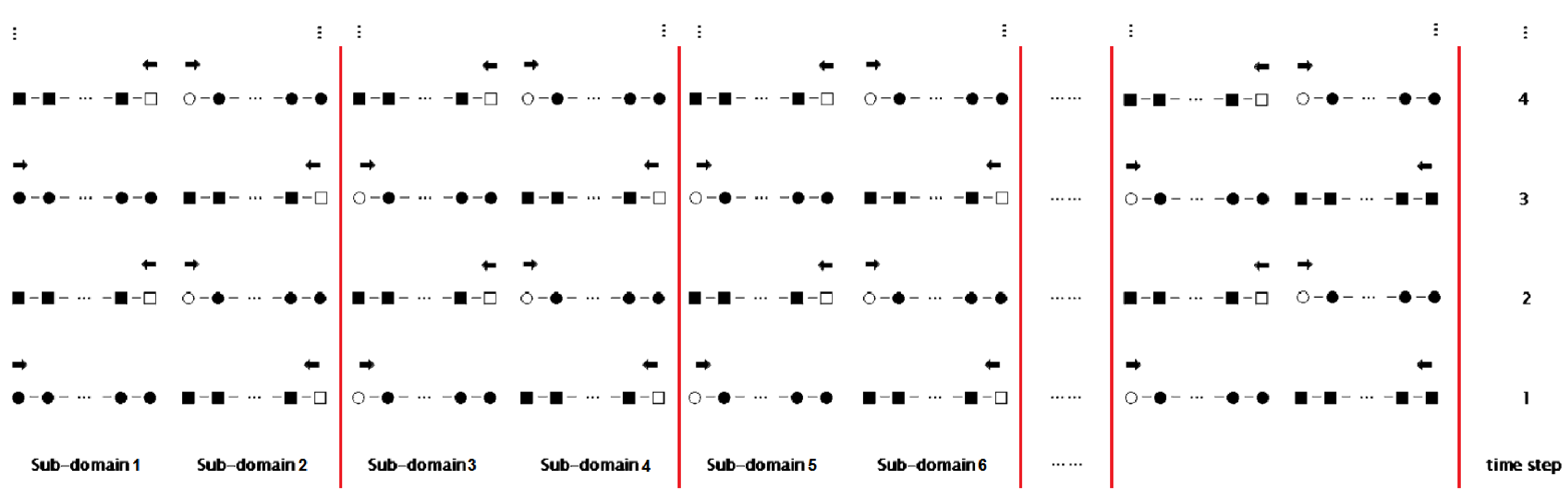}
\caption{Each thread deals with a distinct SDC, determined between red lines}
\label{par_tav_col}
\end{figure}

\begin{figure}
\centering
\includegraphics[scale=1.0]{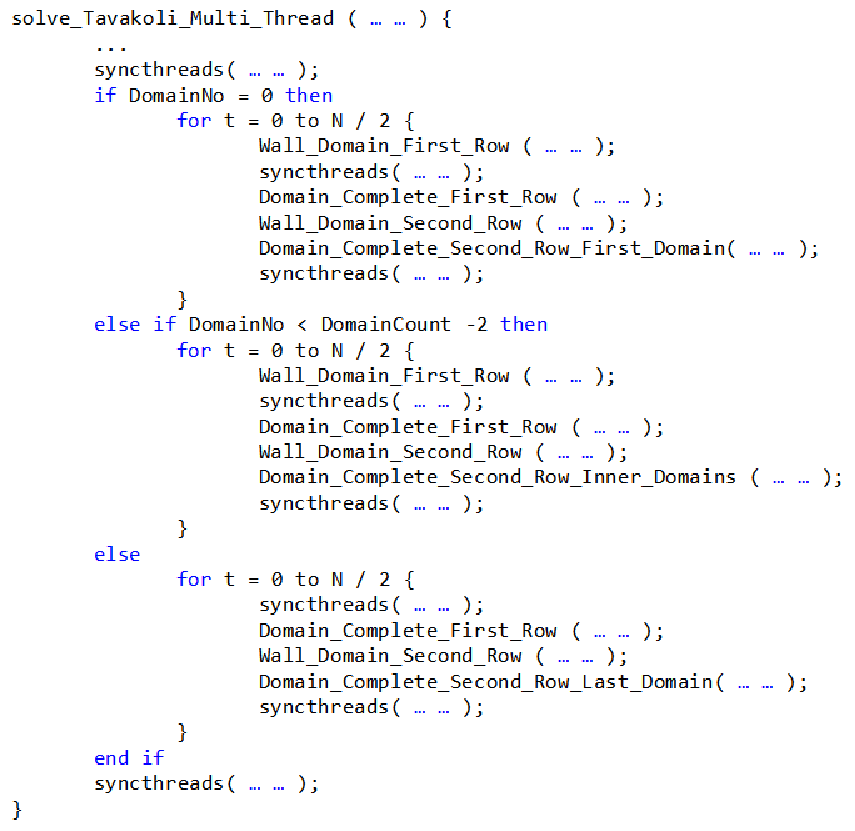}
\caption{The pseudo-code for multi-thread implementation of the parallel scheme}
\label{mtpu}
\end{figure}

In this section, we present a schematic view of multi-thread CPU and GPU implementation of the algorithm and deal with some parallel techniques, for this aim.
As we pointed, before, we chose Tavakoli's\cite{Davam.2007} scheme because of its 4 features:
parallel platform, unconditionally stability, full explicitness and the capacity of extension to higher order dimensions.
As the writer did not find any CPU-implementation report of the diffusion equation solution which satisfies these features, we decided to implement the solution in both CPU and GPU platforms.
In spite of several researches which compare the GPU and single-thread CPU results, we implemented the scheme, one time as a single-thread and the second time as a multi-thread CPU process, to provide much-real comparison results.
As the multi-thread implementation requires some extra techniques more than what we explained in description of the figure \ref{par_tav} (e.g. synchronization) we assigned the rest of this section to explain those techniques.\\
As we mentioned, we are going to present a schematic view of multi-thread CPU and GPU implementation of the algorithm.
As both of the CPU and GPU multi-thread codes follow the same logic, we present a single pseudo-code for both.
The figure \ref{mtpu} shows the implementation pseudo code.
Each GPU or CPU thread should run a copy of the code and deal with a distinct finite difference scheme "Sub-Domain Couple" (SDC), separately.
The SDCs are illustrated, between red lines, in figure \ref{par_tav_col}.\\
We start by introducing the variables.
$N$ is the number of vertical nodes in each column (it is supposed to be an even number);
$DomainCount$ illustrates the number of the total SDCs and $DomainNo$ represents the number of the SDC which the thread is acting on (all of the numberings start from 0).
As it was noted, each independent thread is responsible for one SDC to compute all of its relative column nodes.
But, as it can be understood from the figure \ref{par_tav}, the SDCs are not completely independent.
They use their neighbors' information in the wall nodes.
So, it is necessary to synchronize them at special points.
Accordingly, we applied the $syncthreads()$ function to handle the needed synchronizations.\\
Although $syncthreads()$ function, does not synchronize real blocks and just does synchronization for one virtual block (VB), in the following section and after introducing the concept of VB, we will explain how uniformly distribution of VBs among the stream processors will result in synchronized execution of the entire blocks, in practice.\\
Now, that is the turn of describing the pseudocode described in the figure \ref{mtpu}.
First of all, initialization related to initial and boundary conditions will be done.
Then, each thread goes to a special \texttt{"}if\texttt{"} phrase, depending on its number.
The first thread goes to the first \texttt{"}if\texttt{"} phrase, the last one goes to the last and the other threads go to the second \texttt{"}if\texttt{"} phrase.
Each \texttt{"}if\texttt{"} phrase block is started by a \texttt{"}for\texttt{"} loop which walks on the time dimension, two steps by two steps;
so, at the beginning of each \texttt{"}for\texttt{"} loop iteration, we would, surely, deal with an odd number row of the finite difference scheme.
Also, because of the code synchronized nature, at each iteration, we are sure that all of the threads have the same time ($t$) values (are at the same row).\\
In the first line of \texttt{"}for\texttt{"} loops, the threads begin to compute the nodes of their SDC two right walls ($\Box$ objects as the SDCs right wall and $\Circle$ objects as the SDC right neighbor left wall).
This work is done by calling the $Wall\_Domain\_First\_Row ()$ function (using (\ref{cc}) and (\ref{dd})), except the last SDC which has, only, one preinitialized wall.
As the figure \ref{AB} shows, for computing the wall nodes, each SDC needs some of the down row nodes of its right neighbor;
that is the reason of synchronizing threads before entering the if phrases.\\
Now, each thread computes all of its SDC inner $\CIRCLE$ and $\blacksquare$ objects, from two sides to the SDC center, by calling the $Domain\_Complete\_First\_Row()$ function (using (\ref{aa}) and (\ref{bb})).
It is important to note that the computation of the leftmost SDC $\CIRCLE$ (using \ref{aa}) needs to the SDC left wall node ($\Circle$);
but, computing this node is on the responsibility of the left neighbor thread (using the $Wall\_Domain\_First\_Row ()$ function) and we are not sure if it is computed so far.
So, for ensuring, it is necessary to synchronize the threads, before calling the $Domain\_Complete\_First\_Row()$ function.\\
After running this function, the threads are at their SDC inner wall.
In the next step, unexceptionally, all of the threads compute the upper row middle walls ($\Circle$ and $\Box$ objects) by calling the $Wall\_Domain\_Second\_Row()$ function (using (\ref{cc}) and (\ref{dd})).
Reviewing the formulas (\ref{cc}) and (\ref{dd}) and looking at the figure (\ref{aabb}) clarifies that to compute the $\Box$ and $\Circle$ nodes, we need only two lower nodes of the SDC left sub-domain and two lower nodes of its right sub-domain which are, altogether, computed by calling the recent two discussed functions (note that the 4 needed nodes, for sure, are computed by the current thread, because the minimum size of each SDC width is 4 nodes, as it is obvious in the figure \ref{par_tav_col}).
So, before calling the $Wall\_Domain\_Second\_Row()$ function, no synchronization is needed.\\
Now, that is the turn of computing the SDC middle $\CIRCLE$ and $\blacksquare$ objects from its center to its two sides, using (\ref{aa}) and (\ref{bb}).
But, this process has 2 exceptions: the first and the last SDC.
In the first one, the left wall and in the last one the right wall must not be computed; because they are preinitialized.
To solve this matter, the three functions \texttt{"}Domain\_Complete\_Second\_Row\_First\_Domain()\texttt{"}, \texttt{"}Domain\_Complete\_Second\_Row\_Inner\_Domains()\texttt{"} and \texttt{"}Domain\_Complete\_Second\_Row\_Last\_Domain()\texttt{"} are employed for implementing this processes.
However, it is not necessary to synchronize the threads, before running this process.
Although the computation of the rightmost and the leftmost nodes need their neighbor wall below node, we are sure that they are computed before, because of the synchronized nature of the code, occurred after the below row \texttt{"}wall computing\texttt{"}.\\
We arrived at the end of \texttt{"}for\texttt{"} loops.
It is necessary to synchronize the threads, before finishing the iteration, because the first function  of the next iteration computes the upper row right wall and needs to this row information of the right thread neighbor.
Finishing a \texttt{"}for\texttt{"} loop iteration means that two rows of the finite difference scheme are, completely, computed.
So, by finishing the \texttt{"}if\texttt{"} phrase (finishing all of the \texttt{"}for\texttt{"} loops iterations), the two top rows of the scheme and, consequently, all of the nodes are computed and the algorithm is finished.
\section{Reports and Discussions}
In this section, we report the running results of the discussed algorithm.
We aim to demonstrate the priority of GPU run-time over the CPU.
But, before running the algorithm, there arises a key question; how should we initialize the needed variables?
\subsection{Initializing Variables}
\begin{figure}
\centering
\includegraphics[scale=0.6]{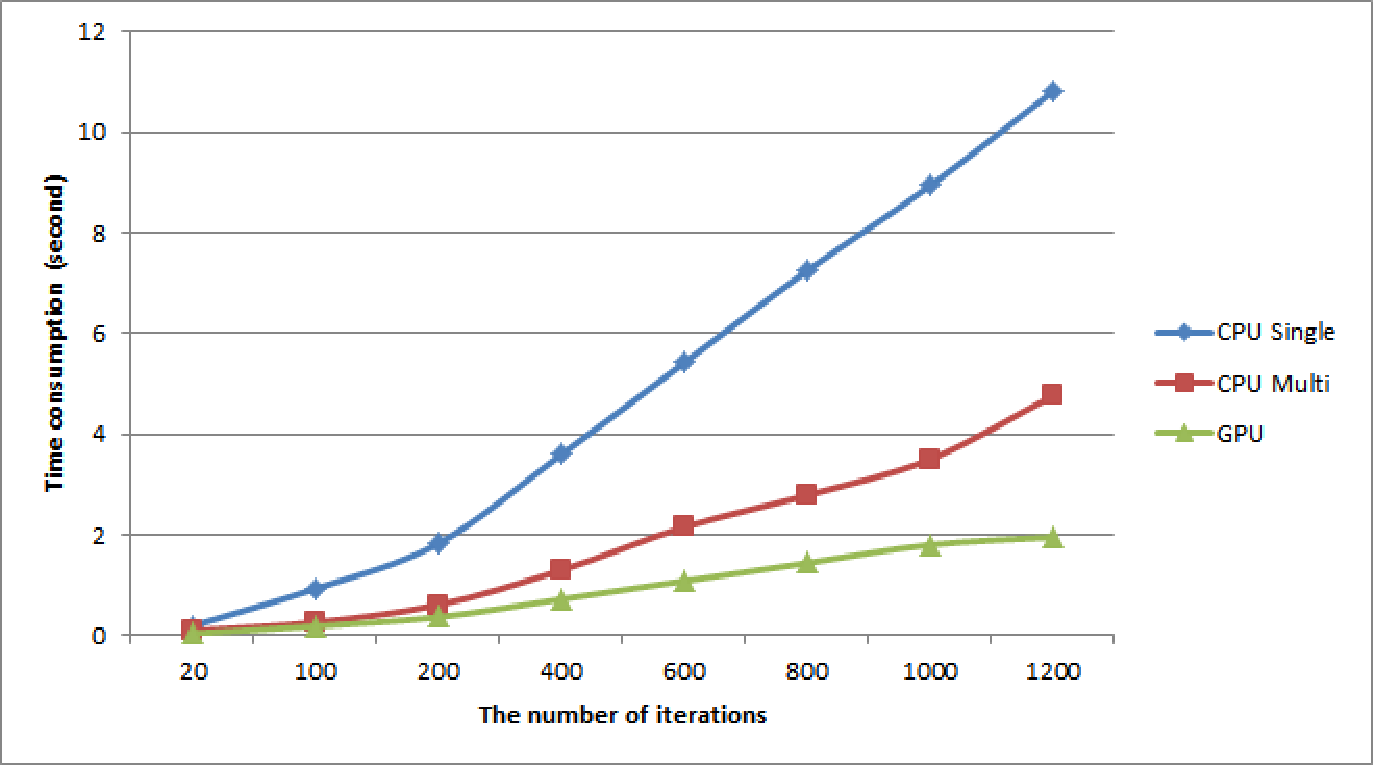}
\caption{The algorithm running time for different threads number}
\label{rsi}
\end{figure}

We should set the number of iterations ($N$, the number of nodes in each column), the number of nodes located in each row ($M$) and the number of threads, for both of the CPU and GPU.\\
Also, to set the number of GPU threads, it is needed to set the number of virtual blocks (VBs) and the number of threads per VB. VBs, as it is understood from the name, are some virtual entities managed by GPU scheduler so that GPU deal with them as real hardware blocks. However, regarding their virtual nature, no processor is assigned to them and each real stream processor is shared equally (and alternatively) with all the VBs which are assigned to that processor.
The CUDA programming model enables the programmer to define an arbitrary number of VBs and each VB can contain from 1 to a maximum of 512 threads
Recall the figure \ref{mtpu} in which the $syncthreads()$ function just synchronizes the threads of one virtual block, since the number of assigned VBs (as we will assign farther (8208)) is a multiple of the number of real blocks (9) and streaming processors (8), GPU scheduler does not cause any time shift in any of the 72 streaming processors and we, experimentally, found that the entire 72 streaming processors work, simultaneously. Note that adopting this approach does not make a limitation in using GPU capabilities. Because, we can add the number of threads, 72 by 72, as much as we need and take advantage of the GPU parallelized nature.\\
As we mentioned, the GPU hardware, used in this work, has 9 real blocks, each of which has 8 stream processors.
We do not set the number of iterations ($N$) and consider it as a variable of the main report plot.
So, by fixing the $M$ variable, we would report various results for different values of finite difference mesh nodes ($M\times N$) and check the effect of the nodes numbers on the run times.\\
As we expressed in the last section, each thread deals with one SDC;
so, setting the thread number is equal to setting the SDC numbers.
Accordingly, the $M$ variable is the product of the thread numbers and the SDC size (the number of SDC row nodes).
We set the SDC size to an arbitrary short power of 2, like 32.
We set it to be short, so that we will be able to deal with a large number of threads and feel the effects of parallel computing in our results, better.
Also, being a power of is, only, for simplifying computations.
Regarding the GPU and CPU capability differences, we set two different thread numbers to each one.
The CPU, used in this work has 8 threads; it means that it has been optimized to run, simultaneously, 8 threads.
So, we set the CPU thread number to 8.
As we explained, setting the GPU thread number is equal to setting the VBs number and the thread per block number.
We set the thread per block variable value to 512, its maximum value, to use whole of the GPU blocks capacity.
But, before setting the VBs number, we performed an investigation on its effect on the GPU running speed.\\
\begin{figure}
\centering
\includegraphics[scale=0.41]{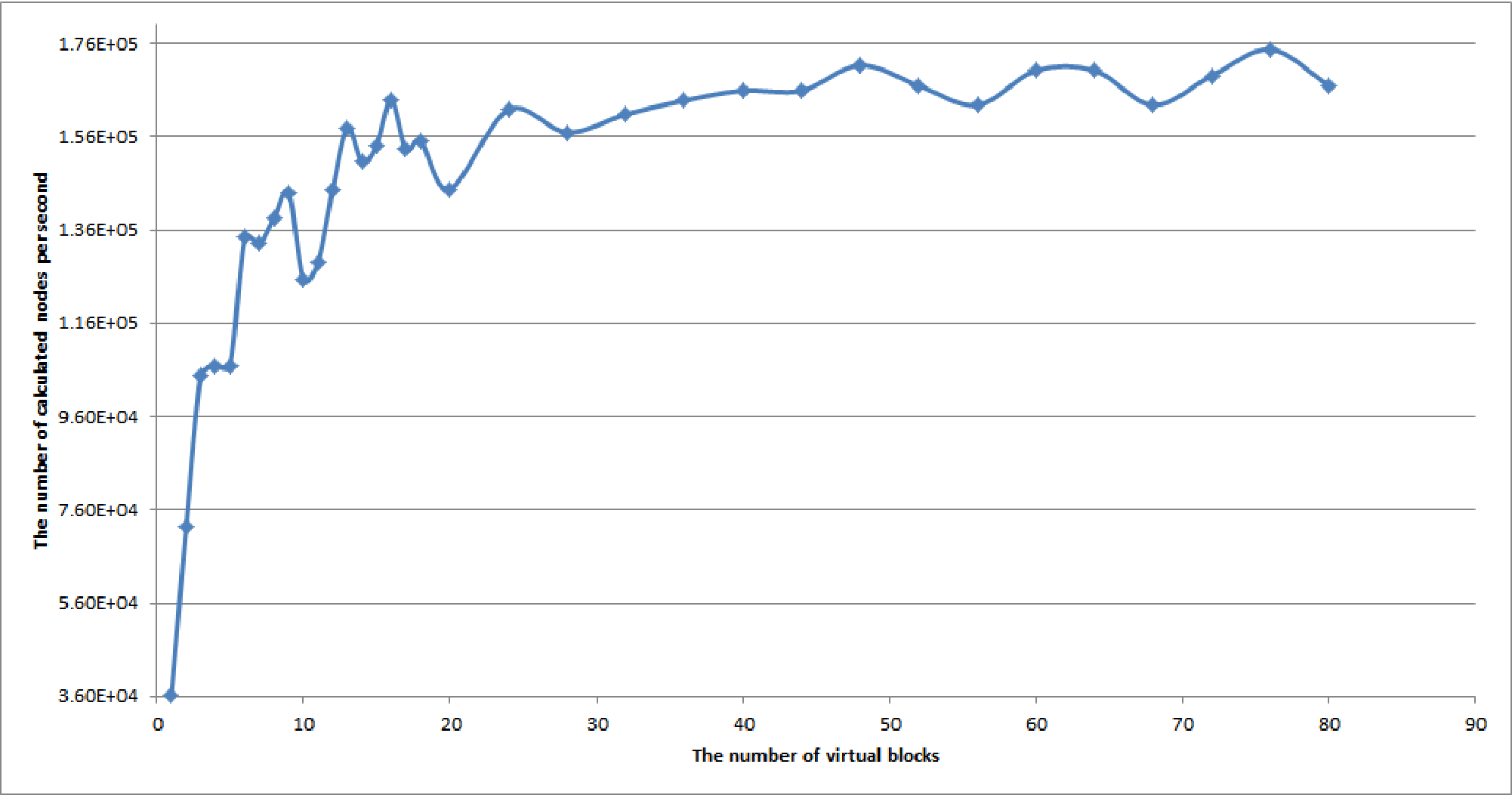}
\caption{The effect of VBs number on the GPU running speed}
\label{bi1880}
\end{figure}
The figure \ref{bi1880} shows the number of calculated finite difference scheme mesh nodes per second (GPU calculation speed) for different VB numbers.
We have run the algorithms for several VB number values and constant nodes numbers (since, VB increment results in $M$ increment and nodes number ($M\times N$) is supposed to be constant, the $N$ number decreases while VB increases.
According to the number of GPU device hardware blocks (9), as the VBs number varies from 1 to 18, we ran the algorithm for each VB number, one by one, to see its effect on the computation speed, more exactly, around 9.
After that we ran the algorithm for multiples of 4, from 20 to 80.
As it was expected, the plot has an increasing manner by increasing the VBs number from 1 to 9, because we are using more GPU hardware blocks.
The interesting event is the plot downfall at the point 9 of the horizontal axis.
It was, completely, expected;
because the pigeonhole principle beside the GPU scheduler says that, in point 10, 1 VB would not be assigned to any hardware block and scheduling algorithm wastages would result in performance decrease.
The plot downfall repeats, again, at the point 18, for the same reason.\\
In general, we can say that increasing the VBs number results in computation speed increment, since, the delivered duties would be broken down into smaller and shorter works;
so, each work needs have little time and, finally, the short time losses in the GPU running process would be revived.
Moreover, about the oscillations it should be said that they are because of the effect of the hardware block numbers on the scheduling management of VBs increment.
Of course, the oscillations decrease by the VBs number increment;
because shrinking the delivered works and reviving the short time GPU losses negate the effect of hardware block numbers, gradually.
Since, we are interested in viewing the results for large $N$ values, it is preferable not to choose large $M$ values;
therefore, we intend to set an almost small VB number which results in an almost high performance.
The figure \ref{bi1880} suggests the VB number to be 16.\\
Therefore, the value of $M$ will be multiplication of number of VBs in threads-per-block, that is $16 \times 512 = 8192$. However, in practice, we change it a little ($8192 + 16 = 8208$) to become divisible into the 72 stream processors. Recall that by distributing the treads, equally, among the 72 stream processors, we made our selves free of synchronizing them. As we tested that in this case, in practice (as we tested), the computations in all of the 72 stream processors are done simultaneously.\\
Now, all of the variables are set and the only remained point is setting a range for $N$.
The experiments showed that our used laptop GPU threads have the ability of computing about 40000 nodes ($M\times N$).
If the nodes number increases from the border, the thread lifetime would finish and the GPU error occurs.
As it is obvious, thread nodes number is the product of its SDC size in the number of iterations ($N$).
Remember that we have set the SDC size to 32;
so the maximum number of iterations which can be set, practically, would be about 1200.
We do not have any limitation for the minimum value and chose 20 for the minimum and ran the algorithm for the set variables and this range.
The results are depicted in the figure \ref{rsi} which show the priority of the GPU run time over the CPU.
It is illustrated that increment in $N$ value does not decrease the GPU superiority.
The GPU times have almost a double excellence over the CPU multi-thread results and a five times excellence over the CPU single thread results.
\subsection{Scalability}
We showed how the cores number and compute capability of the physical processors of a CPU or GPU. The proposed algorithm does not have any limitation, related to any special type of GPUs and by increasing the number of cores, the volume of parallel processes will increase which results in more and more speedup and it is completely independent from other restrictions of a processor, like the copy speed from memory.\\
But letting the GPU to be fixed we cannot scale up the domain of probelm, which is a bottleneck. The entire number of nodes are limited. Therefore if we would like to keep the meshing resolution (needed for results accuracy), we cannot extend the domain in both horizontal and vertical axes; because it results in increment of the number of total nodes whereas our experiments show that GPU cannot suffer this increment.
\section{Conclusion}
In this work, GPU-based accelerating approaches have been applied to stable group explicit finite difference method for solution of one-dimensional diffusion equation.
The method is introduced by \cite{Davam.2007} and has 4 chief features: parallelism platform, unconditional stability, full explicitness (no need to solution of nonlinear algebraic system), capacity of extension to higher dimensions (rather than one).
By using the acceleration, proposed in this work, it would have a fifth feature that is "high speed implementation".
To accomplish the main objective (to evaluate the impact of using GPU in this method), first of all, we discussed on some parallel implementation techniques of the method, like finding the optimized location for synchronizations.
Then, before running the algorithm, we investigated the method of setting the values of the algorithm variables.
Especially, the GPU virtual blocks (VBs) number was set after experimenting the effect of different VB numbers on the GPU speed and, also, by an analytical discussion on the experiment results.
After the value determination of the variables, comparisons between the GPU and the CPU (single and multi-thread) approaches were provided by an ordinary Laptop computer.
For this purpose the CPU multi-thread and single-thread programming codes was implemented, in addition to the GPU code.
The following remarks should be noted as the results:
Using the ordinary laptop computer processors and showing the priority of the GPU results to the CPU results would encourage whole of the researchers and practitioners' society (not especially those who access the HPC labs) to apply GPU on their scientific computing problems.
The GPU results exhibited a speedup about 5 times when compared to a single-thread CPU results and 2 times, when compared to a multi-thread CPU results. In both of the experiments we utilized the entire facilities by which we could keep on writing the simple codes which are familiar with non-computer-expert researchers of this field.
The price of a GT-335M GPU, used in this work (about U\$250) is very cheaper than adding cores in the computer architecture.
The technique has high scalability capacity and there is no limitation for using GPUs with more numerous cores. But, the bottleneck is each GPU has a limitation in total number of nodes which means that we can not increase the problem domain and the grid resolution together, after reaching the maximum number of nodes. 
Finally, it can be stated that this research, like the similar previous studies, demonstrates how easy, cheap and efficient the use of GPU is for solving the scientific computing problems.
\section{Acknowledgment}
We thank Dr.Parviz Davami and, also, Dr.Cyrus Zamani from Iranian Razi Metallurgical Research Center (RMRC), because of their useful consultations.
\bibliographystyle{elsart}
\bibliography{ref}
\end{document}